\documentclass[
final
, nomarks
]{dmtcs-episciences}

\usepackage{mathrsfs}
\usepackage{latexsym,amsmath,amssymb,amsfonts}
\usepackage{color,colordvi,amscd,indentfirst,graphics}
\allowdisplaybreaks
\usepackage[utf8]{inputenc}
\usepackage{subfigure}

\newtheorem{theorem}{Theorem}[section]

\newtheorem{definition}[theorem]{Definition}

\newtheorem{claim}[theorem]{Claim}
\newtheorem{corollary}[theorem]{Corollary}

\def\qed{\hfill $\square$ \nopagebreak}
\def\pf{\noindent {\bf Proof.} }

\usepackage{graphicx}

\author[P. H. Ha and N. G. Hien]{Pham Hoang Ha
  \and Nguyen Gia Hien}
\title[Spanning trees of claw-free graphs with few leaves and branch vertices]{Spanning trees of claw-free graphs with few leaves and branch vertices}
\affiliation{
  Department of Mathematics,  Hanoi National University of Education, Hanoi, Vietnam}
\begin{document}
\publicationdata{vol. 27:3}{2025}{25}{10.46298/dmtcs.15247}{2025-02-17; 2025-02-17; 2025-11-05}{2025-11-08}
\maketitle
\vskip 0.1cm
\begin{abstract}	
 Let $T$ be a tree. A vertex of degree one is a \emph{leaf} of $T$ and a vertex of degree at least three is a \emph{branch vertex} of $T$. A graph is said to be claw-free if it does not contain $K_{1,3}$ as an induced subgraph. In this paper, we study the spanning trees with a bounded number of leaves and branch vertices of claw-free graphs. Applying the main results, we also give some improvements of previous results on the spanning trees with few branch vertices for the case of claw-free graphs.
\end{abstract}
\section{Introduction}
\label{sec:in}
Let $G$ be a finite, simple graph with no loops. The set of vertices and the set of edges of $G$ are denoted by $V(G)$ and $E(G)$, respectively. For each vertex $v$ of $V(G)$, we denote the set of vertices which are adjacent to $v$ in $G$ by $N_G(v)$ and the degree of $v$ in $G$ by $\deg_G(v)$. We define $G-uv$ and $G+uv$ to be the graphs obtained by subtracting and adding the edge $uv$ to $G$, respectively. For every subset $X$ of $V(G)$, the subgraph of $G$ induced by $X$ is denoted by $G[X]$. A graph is called $K_{1,r}$-free if it does not have $K_{1,r}$ as an induced subgraph. A $K_{1,3}$-free graph is also called a claw-free graph.

For a graph $G$, $X\subset V(G)$ is an independent set of $G$ if no two vertices of $X$ are adjacent in $G$. We denote the largest size of independent sets of 
$G$ by $\alpha(G)$. For $k\geq 1$, we define
$$\sigma_k(G)=\left\{\begin{array}{ll}+\infty & \;\mbox{ if } \alpha(G)<k,\\ \min\{\sum_{i=1}^{k}\deg_G(v_i)|\{v_1,\hdots,v_k\} \mbox{ 
		is an independent set of G}\} & \;\mbox{ if } \alpha(G)\geq k.\end{array}\right.$$

Let $T$ be a spanning tree of $G$. A vertex is called a leaf of $T$ if it has degree one in $T$. A vertex is called a branch vertex of $T$ if it has degree strictly greater than two in $T$. The set of leaves and the set of branch vertices of $T$ are denoted by $L(T)$ and $B(T)$, respectively. For each positive integer $k,$ let $B_k(T)$ ($B_{\leq k}(T)$) be the set of branch vertices in $T$ with degree $k$ (at most $k$, respectively).

There are many conditions for a graph $G$ to have a spanning tree $T$ with a bounded number of leaves or branch vertices. We refer the readers to \cite{BT98}, \cite{FKKLR}, \cite{H1},  \cite{Wi} for examples. 

For claw-free graphs, Gargano et al. \cite{GHHSV} proved the following theorem.
\begin{theorem} [Gargano et al. \cite{GHHSV}]
	Let $k\geq 0$ be an integer and let $G$ be a connected claw-free graph. If $\sigma_{k+3}(G)\geq |G|-k-2$, then there exists a spanning tree $T$ of $G$ such that $|B(T)|\leq k$.
\end{theorem}
In 2020, Gould and Shull \cite{GS} proved a conjecture on the spanning tree of a claw-free graph with few branch vertices proposed by Matsuda et al. \cite{MOY}.
\begin{theorem} [Gould and Shull \cite{GS}]
	Let $k\geq 0$ be an integer and let $G$ be a connected claw-free graph. If $\sigma_{2k+3}(G)\geq |G|-2$, then there exists a spanning tree $T$ of $G$ such that $|B(T)|\leq k$.
\end{theorem}
Moreover, many researchers studied the case of $K_{1,r}$-free graphs ($r\geq 4$) , see   \cite{CCHZ}, \cite{CCH14}, \cite{CHH},\cite{Hanh}, \cite{HS}, \cite{Ky09}, \cite{Ky11} for examples.\\
Regarding the conditions for a graph to have a bounded number of leaves and branch vertices, Nikoghosyan \cite{N}, Saito and Sano \cite{SS} independently proved the following.
\begin{theorem} [Nikoghosyan \cite{N}, Saito and Sano \cite{SS}] \label{thm}
	Let $k\geq 2$ be an integer. If a connected graph $G$ satisfies $\sigma_2(G)\geq |G|-k+1$, then there exists a spanning tree $T$ of $G$ such that $|L(T)|+|B(T)|\leq k+1$.
\end{theorem}
In 2019, Maezawa et al. \cite{MMM} gave an improvement of the above result. They proved the following theorem.
\begin{theorem} [Maezawa et al. \cite{MMM}]
	Let $k\geq 2$ be an integer and let $G$ be a connected graph. Suppose that $G$ satisfies $\max\{\deg_G(x),\deg_G(y)\}\geq\dfrac{|G|-k+1}{2}$ for every two vertices $x,y$ such that $xy\notin E(G)$, then there exists a spanning tree $T$ of $G$ such that $|L(T)|+|B(T)|\leq k+1$.
\end{theorem}
For the case of claw-free graphs, Hanh \cite{Hanh1} proved the following theorem. 
\begin{theorem} [Hanh \cite{Hanh1}]
	Let $G$ be a connected claw-free graph. If $\sigma_{5}(G)\geq |G|-2$, then $G$ has a spanning tree with at most $5$ leaves and branch vertices.
\end{theorem}
In the case of $K_{1,4}$-free, Ha \cite{H2} stated the following result.
\begin{theorem} [Ha \cite{H2}]
	Let $G$ be a connected $K_{1,4}$-free graph and $k,m$ be two non-negative integers with $m\leq k+1$. If $\sigma_{m+2}(G)\geq |G|-k$, then there exists a spanning tree $T$ of $G$ such that $|L(T)|+|B(T)|\leq k+m+2$.
\end{theorem}
In the case of $K_{1,5}$-free graphs, two results were introduced as the followings.
\begin{theorem} [Ha and Trang \cite{HT}]
	Let $G$ be a connected $K_{1,5}$-free graph. If $\sigma_4(G)\geq |G|-1$, then there exists a spanning tree $T$ of $G$ such that $|L(T)|+|B(T)|\leq 5$.
\end{theorem}
\begin{theorem} [Diep et al. \cite{DHT}]
	Let $G$ be a connected $K_{1,5}$-free graph. If
	$\sigma_5(G)\geq |G|-2$, then $G$ contains a spanning tree with $|L(T)|+|B(T)|\leq 7$.
\end{theorem}
In this paper, we continue to study some sufficient conditions for a connected claw-free graph
to have a spanning tree with few leaves and branch vertices. The main purpose of this paper is to prove the following theorem.

\begin{theorem}\label{thm-main}
	Let $m,n$ be two positive integers ($n\geq 2$). Let $G$ be a connected claw-free graph. If $\sigma_{m+1}(G)\geq |G|-n+m-1$ and $m \leq \lceil \dfrac{2n}{3}\rceil$, then $G$ has a spanning tree with at most $n$ leaves and branch vertices. Here, the notation $\lceil r\rceil$ stands for the smallest integer not less than the real number $r.$
\end{theorem}

It is easy to see that we directly gain the above result of Hanh \cite{Hanh1} by Theorem \ref{thm-main} with $m=4$ and $n=5$.

Using Theorem \ref{thm-main} with $m=1$ and $n=k+1$ for a positive integer $k,$ then we have the following corollary.
\begin{corollary}\label{Co0}
	Let $k$ be a positive integer and let $G$ be a connected claw-free graph. If $\sigma_{2}(G)\geq |G|-k-1$, then $G$ has a spanning tree with at most $k+1$ leaves and branch vertices.
\end{corollary}
This is an improvement of Theorem \ref{thm} in the case of claw-free graphs.

On the other hand, since $|L(T)|\geq |B(T)|+2$ for each tree $T$, we obtain that if a tree $T$ has at most $2k+3$ leaves and branch vertices, then $|B(T)|\leq k.$ By motivating this fact, we give some sufficient conditions for a claw-free graph to have a spanning tree with few branch vertices.

Let $k$ be an arbitrary positive integer and $n=2k+3.$ Using the same technique of proof of Theorem \ref{thm-main}, we gain the following result.
\begin{theorem}\label{Co1}
	Let $k, m$ be two positive integers such that $k+3\leq m\leq k+\dfrac{k-1}{3}+3$ and let $G$ be a connected claw-free graph. If $\sigma_{m+1}(G)\geq |G|-2k+m-4$, then $G$ has a spanning tree with at most $2k+3$ leaves and branch vertices.
\end{theorem}
In Theorem \ref{Co1}, consider the case $m=k+3,$ we obtain a stronger result of Gargano et al. \cite{GHHSV}.
\begin{corollary}\label{Co2}
	Let $k$ be a positive integer and let $G$ be a connected claw-free graph. If $\sigma_{k+4}(G)\geq |G|-k-1$, then $G$ has a spanning tree with at most $k$ branch vertices.
\end{corollary}

\section{Definitions and Notations}
In this section, we recall some definitions in \cite{GS} which are needed for the proof of main results. We refer to~\cite{Di05} for terminology and notation not defined here.
\begin{definition}[{\cite{GS}}]
	Let $T$ be a tree and let $e$ be an edge of $T$. For any two vertices of $T$, say $u$ and $v$, are joined by a unique path, denoted by $P_T[u,v]$.  We also denote $u_v=V(P_T[u,v])\cap N_T(u)$ and $e_v$ as the vertex incident to $e$ in the direction toward $v$.
\end{definition}

\begin{definition}[{\cite{GS}}]
	Let $T$ be a spanning tree of a graph $G$ and let $v\in V(G)$ and $e\in E(T)$. Denote $g(e,  v)$ as the vertex incident to $e$ farthest away from $v$ in $T$. We say $v$ is an {\it oblique neighbor} of $e$ { with respect to $T$} if $vg(e,  v)\in E(G)$. Let $X\subseteq V(G)$. The edge $e$ has an {\it oblique neighbor in the set} $X$ if there exists a vertex of $X$ which is an  oblique neighbor of $e$  with respect to $T$.
\end{definition}

\begin{definition}[{\cite{GS}}]
	Let $T$ be a spanning tree of a graph $G$. Two vertices are {\it pseudoadjacent} with respect to $T$ if there is some $e\in E(T)$ which has them both as oblique neighbors. Similarly, a vertex set is {\it pseudoindependent} with respect to $T$ if no two vertices in the set are pseudoadjacent with respect to $T$.
\end{definition}

\begin{definition}[{\cite{GS}}]
	Let $T$ be a spanning tree of a graph $G$ with $|B(T)|>0$ and let $r\in B(T)$ be a root of $T$. Then, each branch vertex $b$ has a distance $d(b,r)$ and a degree $\deg_T(b)$. We define a sequence, denoted by $(T,r)$, on the set $B(T)$, which contains the distance-degree pairs of all vertices of $B(T)$ to $r$ in lexicographically increasing order. That is, shortest distance first, and smallest degree first given the same distance.
\end{definition}

\begin{definition}[{\cite{GS}}]  Given two sequences $(T_1,r_1)$ and $(T_2,r_2)$. We define $(T_1,r_1)$ to be smaller than $(T_2,r_2)$ if the distance-degree pair of $(T_1,r_1)$ is smaller than that of $(T_2,r_2)$ at the first different entry between two sequences.
\end{definition}

\section{Proof of Theorem  \ref{thm-main}}


We prove the theorem
by contradiction. Suppose that $G$ has no spanning trees with at most total $n$ leaves and branch vertices. Then, for all spanning trees $T$ of $G$, we have $|L(T)|+|B(T)|\geq n+1$. If $|B(T)|=0$, then $|L(T)|=2$. So $|L(T)|+|B(T)|=2 < n+1$. This is a contradiction. Hence, $|B(T)|\geq 1$. Choose a spanning tree $T$ of $G$ such that:\\
(C1) $|B(T)|$ is as small as possible.\\

We consider two cases as follows.

{\bf Case 1.} $|B_3(T)|=0$ for all spanning trees $T$ satisfying the condition (C1).

In this case, we choose a spanning tree $T$ of $G$ such that:\\
(C2) $|L(T)|$ is as small as possible, subject to (C1).\\
We have $$|L(T)|=2+\sum_{b\in B(T)}(\deg_T(b)-2)\geq 2+2|B(T)|.$$
So $$3|L(T)|\geq 2+2|L(T)|+2|B(T)|\geq 2+2(n+1)=2n+4.$$
Hence, $|L(T)|\geq\dfrac{2n+4}{3}$. Since $|L(T)|$ is an integer and the assumptions of Theorem \ref{thm-main}, we conclude that $|L(T)|\geq m+1$.\\

Let $r\in B(T)$ be the root of $T$.

\begin{claim}\label{claim2}
	For $b\in B(T)$, if $b_1b_2\in E(G)$ and $b_1,b_2$ are children of $b$, then $bb_1$ and $bb_2$ have no oblique neighbors in $L(T)$.
\end{claim}
\pf
Suppose the assertion of the claim is false. Then, there exists $z\in L(T)$ such that $z$ is pseudoadjacent to $bb_1$. If $b_1\in V(P_T[b,z])$, then $T'=T-\{bb_1,bb_2\}+\{bz,b_1b_2\}$ violates the assumption of Case 1 if $b\in B_4(T)$ or condition (C2) otherwise. If $b\in V(P_T[b_1,z])$, then $T'=T-\{bb_1\}+\{zb_1\}$ violates the assumption of Case 1 if $b\in B_4(T)$ or condition (C2) otherwise. The case for $bb_2$ is done by symmetry. This completes the proof of Claim \ref{claim2}.
\qed

\begin{claim}\label{claim3}
	$L(T)$ is an independent set.
\end{claim}
\pf
Suppose that there are two leaves $u,v$ of $T$ such that $uv\in E(G)$. Let $t$ be the nearest branch
vertice of $u.$ Then $T'=T-\{tt_u\}+\{uv\}$ violates either the assumption of Case 1 if $t\in B_4(T)$ or the condition (C2) for otherwise. Therefore, Claim \ref{claim3} is proved.
\qed

\begin{claim}\label{claim4}
	$L(T)$ is a pseudoindependent set.
\end{claim}
\pf
Suppose to the contrary that there exists $\{u,v\}\subset L(T)$ and an edge $e$ of $T$ such that $ug(e,u)\in E(G)$ and $vg(e,v)\in E(G)$. \\
If $g(e,u)=g(e,v)=a$. Denote $\{t\}=V(P_T[u,a])\cap V(P_T[a,v])\cap V(P_T[v,u])$. Since $G[a,e_t,u,v]$ is not a claw and $uv\notin E(G)$, we obtain either $ue_t\in E(G)$ or $ve_t\in E(G)$. Without loss of generality, we may assume that $ue_t\in E(G)$. Then the spanning tree $T'=T-\{e,tt_u\}+\{ue_t,va\}$ contradicts either the assumption of Case 1 if $t\in B_4(T)$ or the condition (C2) if not.\\

\noindent If $g(e,u)=e_v$ and $g(e,v)=e_u$, then $e\in P_T[u,v]$.
 Denote $\{s\}=V(P_T[u,r])\cap V(P_T[r,v])\cap V(P_T[v,u])$. Without loss of generality, we may assume that $s\in V(P_T[e_u,u])$. Then $T'=T-\{e,ss_u\}+\{ue_v,ve_u\}$ violates the assumption of Case 1 if $s\in B_4(T)$ or the condition (C2) otherwise. We conclude that $L(T)$ is a pseudoindependent set.
\qed\\

Let $Q$ be a subset of $L(T)$ such that $|Q|=m+1$. Let $s$ be the number of edges in $T$ which have no oblique neighbors in $Q$. 

By Claim \ref{claim3}, we deduce that for each leaf $l\in L(T)\setminus Q$, $ll_r$ has no oblique neighbors in $Q$. Let $A$ be the set of all such edges $ll_r$ in $T$.

On the other hand, for each $b\in B(T),$ $\deg_T(b)\geq 4$ and the fact that $G$ is claw-free, we obtain that there exist two children $b_1,b_2$ of $b$ such that $b_1b_2\in E(G).$ Then $bb_1$ and $bb_2$ have no oblique neighbors in $Q$ from the fact of Claim \ref{claim2}. We denote by $B$ the set of all such edges $bb_1, bb_2$ in $T$.

We will prove that $A$ and $B$ are disjoint. Suppose that there exists an edge $e$ of $T$ such that $e\in A\cap B$. We deduce that $e_r\in B(T)$, $g(e,r)\in L(T)$ and $g(e,r)$ is adjacent to a child $z\neq g(e,r)$ of $e_r$. Then $T'=T-\{e_rz\}+\{zg(e,r)\}$ violates the assumption of Case 1 if $e_r\in B_4(T)$ or condition (C2) otherwise. We conclude that $A$ and $B$ are disjoint sets.

Therefore, we have
\begin{align*}
	s&\geq |A|+|B|\geq|L(T)|-|Q|+2|B(T)\setminus\{r\}|+2=|L(T)|-(m+1)+2|B(T)|\\
	&\geq n+1-m.
\end{align*}

On the other hand, for any $x, y\in V(T)$, we have $xy\in E(G)$  if and only if  $x$ is an oblique neighbor of $yy_x$.  Therefore, the number of edges of $T$ with $x$ as an oblique neighbor equals the degree of $x$ in $G$. Therefore, combining with Claim \ref{claim4}, we obtain that
$$\sigma_{m+1}(G)\leq\sum_{t\in Q}\deg_G(t)\leq |E(T)|-s\leq (|G|-1)-(n-m+1)=|G|-n+m-2.$$
This contradicts with the assumption of Theorem \ref{thm-main}. \\
\vskip 0.5cm
{\bf Case 2.} There exists at least one spanning tree $T$ of $G$ such that $|B_3(T)|>0$.

Let $r\in B_3(T)$ be a root of $T.$

In this case, in all spanning trees satisfying the condition (C1) and having at least one branch vertex of degree $3$, we choose a spanning tree $T$ with the root $r$ such that:\\
(C3) $\sum_{v\in B_{\geq 5}(T)} (\deg_T(v)-4)$ is as small as possible.\\
(C4) $(T,r)$ is lexicographically as small as possible, subject to (C3).\\
Thus $\deg_T(r)=3$.\\
We have $$|L(T)|=2+\sum_{b\in B(T)}(\deg_T(b)-2)\geq 2+|B_3(T)|+2|B_{\geq 4}(T)|.$$
So $$
	|L(T)|+|B_3(T)| \geq 2+2|B_3(T)|+2|B_{\geq 4}(T)|
	=2(1+|B(T)|)\geq 2(1+n+1-|L(T)|).
$$
This implies  $$ 3|L(T)|+3|B_3(T)|\geq 2n+4+2|B_3(T)|\geq 2n+6.$$
Therefore, we obtain $$ |L(T)|+|B_3(T)| \geq\frac{2n+6}{3}\geq m+2.$$
Since $|L(T)|+|B_3(T)|$ is an integer and the assumptions of Theorem \ref{thm-main}, we obtain $|L(T)|+|B_3(T)|\geq m+2$.
Let $H=L(T)\cup B_3(T)\setminus\{r\}$. Then $|H|\geq m+1.$ We now have the following claims.

\begin{claim}\label{claim5}
	If $u\in B(T)\setminus\{r\}$ and $a$ is a child of $u$, then $a$ is adjacent to at least one neighbor of $u$.
\end{claim}
\pf
Assume that $ab\notin E(G)$ for all $b\in N_T(u)\setminus\{a\}$. Then let $c$ be a child of $u$ which is different from $a$. Since $G[u,a,c,u_r]$ is claw-free and $au_r,ac\notin E(G)$, we have $cu_r\in E(G)$. This concludes that $bu_r\in E(G)$ for every $b\in N_T(u)\setminus\{a,u_r\}.$ Then the spanning tree $T'=T-\{ub|b\in N_T(u)\setminus\{a,u_r\}\}+\{u_rb|b\in N_T(u)\setminus\{a,u_r\}\}$ violates either the condition (C1) if $u_r\in B(T)$ or the condition (C4) if $u_r\notin B(T)$. So the claim holds.
\qed

\begin{claim}\label{claim6}
	If $u\in B_3(T)\setminus\{r\}$ and $a,b$ are two children of $u$, then $ab\in E(G)$.
\end{claim}
\pf
Suppose for a contradiction that $ab\notin E(G)$, since $G[u,a,b,u_r]$ is not a claw and $ab\notin E(G)$, we obtain either $au_r\in E(G)$ or $bu_r\in E(G)$. Without loss of generality, we may assume that $au_r\in E(G)$. Then the spanning tree $T'=T-\{au\}+\{au_r\}$ contradicts either the condition (C1) if $u_r\in B(T)$ or the condition (C4) if $u_r\notin B(T)$. This completes the proof of Claim \ref{claim6}.
\qed

\begin{claim}\label{claim7}
	$H$ is an independent set.
\end{claim}
\pf
Suppose this is false. Then there exists $\{u,v\}\subset H$ such that $uv\in E(G)$. If $u\in V(P_T[r,v])$, then $\deg_T(u)=3$. We denote $\{u^*\}=N_T(u)\setminus\{u_r,u_v\}$. By Claim \ref{claim6}, we obtain that $u^*u_v\in E(G)$. Then $T'=T-\{uu^*,uu_v\}+\{u^*u_v,uv\}$ violates the condition (C1). The case $v\in V(P_T[r,u])$ is done by symmetry. Otherwise, we have $\{w\}=V(P_T[r,u])\cap V(P_T[u,v])\cap V(P_T[v,r])\not\subset \{u,v\}.$ Consider the spanning tree $T'=T-\{ww_u\}+\{uv\}.$ This contradicts either the condition (C1) if $w\in B_3(T)$, or the condition (C3) if $w\in B_{\geq 5}(T)$ or the condition (C4) if $w\in B_4(T)$. So we conclude that $H$ is an independent set.
\qed

\begin{claim}\label{claim8}
	$H$ is a pseudoindependent set.
\end{claim}
\pf
Assume that there exists $\{u,v\}\subset H$ and an edge $e\in E(T)$ such that $ug(e,u)\in E(G)$ and $vg(e,v)\in E(G)$. Let $x$ be a leaf or branch vertex which is nearest to $e$ in the direction away from $r$. Denote $\{w\}=V(P_T[u,v])\cap V(P_T[v,r])\cap V(P_T[r,u])$.

If $e\in P_T[u,v]$, then $g(e,u)\neq g(e,v)$. Now if $u\in V(P_T[r,v])$, then $\deg_T(u)=3$ and we denote $\{u^*\}=N_T(u)\setminus\{u_v,u_r\}$. By Claim \ref{claim6}, we may obtain that $u^*u_v\in E(G)$. Then $T'=T-\{uu_v,uu^*,e\}+\{ug(e,u),vg(e,v),u^*u_v\}$ violates the condition (C1). The case $v\in V(P_T[r,u])$ is done by symmetry. If $w\notin\{u,v\}$, we consider the spanning tree
$$T'=\left\{\begin{array}{ll}T-\{e,ww_u\}+\{ug(e,u),vg(e,v)\}, & \;\mbox{ if } e\neq ww_u,\\ T-\{e\}+\{vw_u\}, & \;\mbox{ if } e=ww_u.\end{array}\right.$$
Then $T'$ violates either the condition (C1) if $w\in B_3(T)$, or the condition (C3) if $w\in B_{\geq 5}(T)$ or the condition (C4) if $w\in B_4(T)$.

If $e\notin P_T[u,v]$, denote $\{p\}=V(P_T[x,u])\cap V(P_T[u,r])\cap V(P_T[r,x])$ and $\{q\}=V(P_T[x,v])$ $\cap V(P_T[v,r])\cap V(P_T[r,x])$. We consider four cases as follows.

Case \ref{claim8}.1: Suppose $r\in V(P_T[w,e_x])$, then $p=r$ or $q=r$. Without loss of generality, we may assume that $p=r$, thus $g(e,u)=g(e,v)=e_x$. Since $G[e_x,e_r,u,v]$ is claw-free, we obtain either $ue_r\in E(G)$ or $ve_r\in E(G)$. 

If $ue_r\in E(G)$, we consider the spanning tree $$T'=\left\{\begin{array}{ll}T-\{e,rr_x\}+\{ue_r,ve_x\}, & \;\mbox{ if } e\neq rr_x,\\ T-\{e\}+\{ve_x\}, & \;\mbox{ if } e=rr_x.\end{array}\right.$$
Then $T'$ violates the condition (C1).

If $ve_r\in E(G)$, we will prove that $e_x\neq x$. Assume that $e_x=x$, then $T'=T-\{rr_u\}+\{ue_x\}$ violates the condition (C1), thus $e_{xx}$ exists. Since $G[e_x,e_{xx},u,v]$ is claw-free, we obtain either $ue_{xx}\in E(G)$ or $ve_{xx}\in E(G)$. If $ve_{xx}\in E(G)$, then $T'=T-\{e_xe_{xx},rr_u\}+\{ue_x,ve_{xx}\}$ violates the condition (C1). If $ue_{xx}\in E(G)$, then $T'=T-\{e_xe_{xx},qq_x\}+\{ve_x,ue_{xx}\}$ violates either the condition (C1) if $q\in B_3(T)$, or the condition (C3) if $q\in B_{\geq 5}(T)$ or the condition (C4) if $q\in B_4(T)$.

Case \ref{claim8}.2: Suppose $e_x\in V(P_T[r,w])$, then $g(e,u)=g(e,v)=e_r$. Since $G[e_r,e_x,u,v]$ is claw-free, we obtain either $ue_x\in E(G)$ or $ve_x\in E(G)$. If $u\in V(P_T[r,v])$, then $T'=T-\{uu_v\}+\{ve_r\}$ violates either the condition (C1) if $e_r\in B(T)$ or the condition (C4) otherwise. The case $v\in V(P_T[r,u])$ is done by symmetry. Now we consider the case $w\notin\{u,v\}$. Without loss of generality, we may assume that $ue_x\in E(G)$. If $w\in B_3(T)$, then $T'=T-\{ww_u\}+\{ue_r\}$ violates the condition (C1) if $e_r\in B(T)$ or condition (C4) if not. If $w\in B_4(T)$, then $T'=T-\{ww_u,ww_v\}+\{ue_r,ve_r\}$ violates the condition (C1) if $e_r\in B(T)$ or condition (C4) if not. If $w\in B_{\geq 5}(T)$, then $T'=T-\{e,ww_u\}+\{ue_x,ve_r\}$ violates the condition (C3).

Case \ref{claim8}.3: Suppose $w\in V(P_T[r,e_r])\setminus\{r\}$ and $e\notin P_T[u,v]$, then $p=w$ or $q=w$. Without loss of generality, we may assume that $p=w$, thus $g(e,u)=g(e,v)=e_x$. We consider two cases as follows.

Subcase 1: $w\in \{u, v\}.$ By symmetry, we also assume that $w=u$. This implies $u\in V(P_T[r,v]).$ We will prove that $e_x\neq x$. Denote $\{u^*\}=N_T(u)\setminus\{u_r,u_x\}$. By Claim \ref{claim6}, we obtain that $u^*u_x\in E(G)$. If $e_x=x$, then $T'=T-\{uu^*,uu_x\}+\{u^*u_x,ux\}$ violates the condition (C1), so $e_{xx}$ exists. Since $G[e_x,e_{xx},u,v]$ is not a claw graph and $uv\notin E(G)$, we may obtain either $ue_{xx}\in E(G)$ or $ve_{xx}\in E(G)$. If $ve_{xx}\in E(G)$, then $T'=T-\{e_xe_{xx},uu_x,uu^*\}+\{u_xu^*,ve_{xx},ue_x\}$ violates the condition (C1). If $ue_{xx}\in E(G)$, we now consider the degree and position of $q$. 

If $q=v$, denote $\{v^*\}=N_T(v)\setminus\{v_x,v_r\}$. By Claim \ref{claim6}, we obtain that $v^*v_x\in E(G)$. Then $T'=T-\{uu^*,uu_x,vv^*,vv_x\}+\{u^*u_x,v^*v_x,ue_x,ve_x\}$ violates the condition (C1).

If $q=u$, since $G[u,u_r,u_x,e_{xx}]$ is claw-free and $u_ru_x\notin E(G)$, we obtain either $u_re_{xx}\in E(G)$ or $u_xe_{xx}\in E(G)$. If $u_re_{xx}\in E(G)$, then $T'=T-\{e_xe_{xx},uu^*\}+\{u_re_{xx},ve_x\}$ violates either the condition (C1) if $u_r\in B(T)$ or the condition (C4) if $u_r\notin B(T)$. If $u_xe_{xx}\in E(G)$, then $T'=T-\{e_xe_{xx},uu_x\}+\{u_xe_{xx},ve_x\}$ violates the condition (C1).

Now if $q\notin \{u,v\}.$ Then $T'=T-\{qq_v,e_xe_{xx}\}+\{ue_{xx},ve_x\}$ violates either the condition (C1) if $q\in B_3(T)$ or the condition (C3) if $q\in B_{\geq 5}(T).$ So we only have to consider the case $q\in B_4(T)$. By Claim \ref{claim5}, $q_v$ is adjacent to a vertex in $N_T(q)\setminus\{q_v\}$. Denote $\{q^*\}=N_T(q)\setminus\{q_x,q_v,q_r\}$. We consider the spanning tree
$$T'=\left\{\begin{array}{ll}T-\{e_xe_{xx},qq_v,qq_r\}+\{q_vq_r,ve_x,ue_{xx}\}, & \;\mbox{ if } q_vq_r\in E(G),\\ T-\{e_xe_{xx},qq_v,qq_x\}+\{q_vq_x,ue_x,ue_{xx}\}, & \;\mbox{ if } q_vq_x\in E(G),\\ T-\{e_xe_{xx},qq_v,qq^*\}+\{q_vq^*,ve_x,ue_{xx}\}, & \;\mbox{ if } q_vq^*\in E(G).\end{array}\right.$$
Then $T'$ violates the condition (C1).

Subcase 2: $w\notin\{u,v\}$. Since $G[e_x,e_r,u,v]$ is claw-free, we obtain either $ue_r\in E(G)$ or $ve_r\in E(G)$. 

If $p=q=w$, without loss of generality, we may assume that $ue_r\in E(G)$. We consider the spanning tree
$$T'=\left\{\begin{array}{ll}T-\{e,ww_x\}+\{ue_r,ve_x\}, & \;\mbox{ if } e\neq ww_x,\\ T-\{e\}+\{ve_x\}, & \;\mbox{ if } e=ww_x.\end{array}\right.$$
Then $T'$ violates the condition (C1) if $w\in B_3(T)$, the condition (C3) if $w\in B_{\geq 5}(T)$ or the condition (C4) if $w\in B_4(T)$.

If $q=v$, denote $\{v^*\}=N_T(v)\setminus\{v_x,v_r\}$. We will prove that $e_x\neq x$. Suppose that $e_x=x$, by Claim \ref{claim6}, we obtain that $v^*v_x\in E(G)$. Then $T'=T-\{vv_x,vv^*\}+\{ve_x,v_xv^*\}$ violates the condition (C1), so $e_x\neq x$. Now $T'=T-\{ww_x,vv_x,vv^*\}+\{v^*v_x,ue_x,ve_x\}$ violates the condition (C1) if $w\in B_3(T)$, the condition (C3) if $w\in B_{\geq 5}(T)$ or the condition (C4) if $w\in B_4(T)$.

If $q\neq w$ and $q\neq v$, we will prove that $e_x\neq x$. Suppose that $e_x=x$, by Claim \ref{claim7}, we deduce that $\deg_T(x)\geq 4$. Then $T'=T-\{qq_v\}+\{ve_x\}$ violates either the condition (C1) if $q\in B_3(T)$ or the condition (C3) if $q\in B_{\geq 5}(T).$ So we only have to consider the case $q\in B_4(T)$. By Claim \ref{claim5}, $q_v$ is adjacent to a vertex in $N_T(q)\setminus\{q_v\}$. Denote $\{q^*\}=N_T(q)\setminus\{q_x,q_v,q_r\}$. We consider the spanning tree
$$T'=\left\{\begin{array}{ll}T-\{qq_v,qq_r\}+\{q_vq_r,ve_x\}, & \;\mbox{ if } q_vq_r\in E(G),\\ T-\{qq_v,qq_x\}+\{q_vq_x,ue_x\}, & \;\mbox{ if } q_vq_x\in E(G),\\ T-\{qq_v,qq^*\}+\{q_vq^*,ve_x\}, & \;\mbox{ if } q_vq^*\in E(G).\end{array}\right.$$
Then $T'$ violates the condition (C1). So $e_x\neq x$, thus $e_{xx}$ exists. 

Since $G[e_x,e_{xx},u,v]$ is claw-free and $uv\notin E(G)$, we may obtain either $ue_{xx}\in E(G)$ or $ve_{xx}\in E(G)$. If $ve_{xx}\in E(G)$, then $T'=T-\{e_xe_{xx},ww_x\}+\{ue_x,ve_{xx}\}$ contradicts either the condition (C1) if $w\in B_3(T)$ or the condition (C3) if $w\in B_{\geq 5}(T)$ or the condition (C4) if $w\in B_4(T)$. If $ue_{xx}\in E(G)$, then $T'=T-\{e_xe_{xx},qq_v\}+\{ve_x,ue_{xx}\}$ violates either the condition (C1) if $q\in B_3(T)$ or the condition (C3) if $q\in B_{\geq 5}(T)$. Now we consider the case $q\in B_4(T)$. By Claim \ref{claim5}, $q_v$ is adjacent to a vertex in $N_T(q)\setminus\{q_v\}$. Denote $\{q^*\}=N_T(q)\setminus\{q_x,q_v,q_r\}$. We consider the spanning tree
$$T'=\left\{\begin{array}{ll}T-\{e_xe_{xx},qq_v,qq_r\}+\{q_vq_r,ve_x,ue_{xx}\}, & \;\mbox{ if } q_vq_r\in E(G),\\ T-\{e_xe_{xx},qq_v,qq_x\}+\{q_vq_x,ue_x,ue_{xx}\}, & \;\mbox{ if } q_vq_x\in E(G) \mbox{ and } u\in B_3(T),\\ T-\{e_xe_{xx},qq_v,qq^*\}+\{q_vq^*,ve_x,ue_{xx}\}, & \;\mbox{ if } q_vq^*\in E(G).\end{array}\right.$$
Then $T'$ violates the condition (C1). If $q_vq_x\in E(G)$ and $u\in L(T)$, then $T'=T-\{e_xe_{xx},qq_v\}+\{ve_x,ue_{xx}\}$ violates the condition (C1) if $q_v\in B_3(T)$, the condition (C3) if $q_v\in B_{\geq 5}(T)$ or the condition (C4) otherwise.

Case \ref{claim8}.4: If no vertex in the set $\{w,r,e_x\}$ is between the other two vertices, then $ue_x\in E(G), ve_x\in E(G)$ and $p=q$. Since $G[e_x,e_r,u,v]$ is claw-free, we obtain either $ue_r\in E(G)$ or $ve_r\in E(G)$. Without loss of generality, we may assume that $ue_r\in E(G)$. We consider the spanning tree
$$T'=\left\{\begin{array}{ll}T-\{e,qq_x\}+\{ue_r,ve_x\}, & \;\mbox{ if } e\neq qq_x,\\ T-\{e\}+\{ue_x\}, & \;\mbox{ if } e=qq_x.\end{array}\right.$$
Then $T'$ violates the condition (C1) if $q\in B_3(T)$, the condition (C3) if $q\in B_{\geq 5}(T)$ or the condition (C4) if $q\in B_4(T)$.

So we conclude that $H$ is a pseudoindependent set.
\qed

\begin{claim}\label{claim9}
	If $r_1,r_2$ are adjacent to $r$ in $T$ and $r_1r_2\in E(G)$, then $rr_1$ and $rr_2$ have no oblique neighbors in $H$.
\end{claim}
\pf
Suppose that there exists $z\in H$ and $z$ is pseudoadjacent to $rr_1$. We consider the spanning tree
$$T'=\left\{\begin{array}{ll}T-\{rr_1, rr_2\}+\{rz,r_1r_2\},
	& \;\mbox{ if } r_1\in V(P_T[r,z]),\\
	T-\{rr_1\}+\{zr_1\}, & \;\mbox{ if } r\in V(P_T[r_1,z]).
\end{array}\right.$$
Then $T'$ violates condition (C1). Hence $rr_1$ has no oblique neighbors in $H$. By symmetry, $rr_2$ has no oblique neighbors in $H$.
\qed\\

Since $|H|\geq m+1$, we choose a subset $M$ of $H$ such that $|M|=m+1.$ 
\begin{claim}\label{claim10}
	For $b\in B(T)\setminus (M\cup\{r\})$, if $b_1b_2\in E(G)$ and $b_1,b_2$ are children of $b$, then $bb_1$ or $bb_2$ has no oblique neighbors in $M$.
\end{claim}
\pf
Suppose to the contrary that there exists some vertex $b\in B(T)\setminus (M\cup\{r\})$ with two children $b_1,b_2$ and $z, t\in M$ such that $b_1b_2\in E(G)$ and $z, t$ are pseudoadjacent to $bb_1, bb_2,$ respectively. We consider two cases as follows.

Case \ref{claim10}.1: Suppose $b_1\in V(P_T[b,z])$. Then $T'=T-\{bb_1,bb_2\}+\{bz,b_1b_2\}$ violates either the condition (C1) if $b\in B_3(T)$, or the condition (C3) if $b\in B_{\geq 5}(T)$ or the condition (C4) if $b\in B_4(T)$.

Case \ref{claim10}.2: Suppose $b\in V(P_T[b_1,z])$. Then $T'=T-\{bb_1\}+\{zb_1\}$ violates either the condition (C1) if $b\in B_3(T)$ or the condition (C3) if $b\in B_{\geq 5}(T)$. Now we consider the case $b\in B_4(T)$ and $b\in V(P_T[b_1,z])\cap V(P_T[b_2,t])$. If $z\neq t$, $b_1\notin V(P_T[t,b])$ and $b_2\notin V(P_T[z,b])$, then $T'=T-\{bb_1,bb_2\}+\{zb_1,tb_2\}$ violates the condition (C1). Now if $z\neq t$ and $b_1\in V(P_T[t,b])$, then $T'=T-\{bb_2\}+\{tb_2\}$ violates the condition (C4). If $z=t$ and $d(b,r)\geq d(z,r)$, then $T'=T-\{bb_1,bb_2\}+\{zb_1,zb_2\}$ violates either the condition (C1) if $z\in B_3(T)$ or the condition (C4) if $z\in L(T)$, as $\deg_{T'}(z)=3<4=\deg_{T}(b)$ and $d(b,r)\geq d(z,r)$. If $z=t$ and $d(b,r)<d(z,r)$, then $T'=T-\{bb_1\}+\{zb_1\}$ violates the condition (C4). This completes the proof of Claim \ref{claim10}.
\qed\\

Now we conclude that:

If $u\in B_3(T)\setminus (M\cup\{r\}) ,$ combining with Claims \ref{claim6} and \ref{claim10}, then there exists at least one child $a$ of $u$ such that the edge $au$ has no oblique neighbors in $M.$ If $u=r$, there exist at least two neighbors $r_1, r_2$ of $r$ such that $rr_1, rr_2$ have no oblique neighbors in $M$ by Claim \ref{claim9}. Let $C$ be the set of all such edges in $T$.

If $b\in B_{\geq 4}(T) .$ Since $G$ is claw-free, $b$ has two children $b_1,b_2$ such that $b_1b_2\in E(G).$ Combining with Claim \ref{claim10}, we obtain that there exists at least one child $c$ of $b$ such that the edge $bc$ has no oblique neighbors in $M$ and $c$ is adjacent to at least one different child of $b$. Let $D$ be the set of all such edges in $T$. 

By Claim \ref{claim7}, we deduce that for each leaf $l\in L(T)\setminus M$, $ll_r$ has no oblique neighbors in $M$. Let $E$ be the set of all such edges in $T$.

We will prove that $C$, $D$ and $E$ are disjoint sets. The fact that $C \cap D = \emptyset$ is trivial, so suppose to the contrary that there exists $e\in E(T), e\in (C\cup D)\cap E$. So, we deduce that $e_r\in B(T), g(e,r)\in L(T)$ and $g(e,r)$ is adjacent to a child $z\neq g(e,r)$ of $e_r$. We consider the tree $T'=T-\{e_rz\}+\{zg(e,r)\}$. This spanning tree violates either the condition (C1) if $e_r\in B_3(T)$, or the condition (C3) if $e_r\in B_{\geq 5}(T)$ or the condition (C4) otherwise. Then we conclude that $C$, $D$ and $E$ are distinct.

Therefore, if we set $h$ to be the number of edges in $T$ which have no oblique neighbors in $M$, then we get
\begin{align*}
	h& \geq |C|+|D|+|E|\\
	&\geq 2+|(B(T)\setminus\{r\})\setminus M)|+|L(T)\setminus M|\\
	&=|B(T)|+1+|L(T)|-|M|=|L(T)|+|B(T)|-m\\
	&\geq n+1-m.
\end{align*}
So, we obtain that
$$\sigma_{m+1}(G)\leq\sum_{t\in M}\deg_G(t)\leq |E(T)|-h\leq |G|-n+m-2.$$
This gives a contradiction with the assumption of Theorem \ref{thm-main}.
The proof of Theorem \ref{thm-main} is completed.



\end{document}